\documentclass{svproc}
\usepackage{url}

\usepackage[textsize=footnotesize backgroundcolor=yellow!70, bordercolor=orange]{todonotes}
\usepackage{hyperref} 
\usepackage{amsmath}
\usepackage{amssymb}
\usepackage{amsfonts}
\usepackage{subcaption}
\usepackage{graphicx,wrapfig,lipsum}
\newcommand{\dt}{\Delta t}

\makeatletter

\newcommand{\ha}{\frac{1}{2}}

\renewcommand{\epsilon}{\varepsilon }

\newfont{\numerikEleven}{ecrm1000}
\newfont{\numerikTen}{cmss10}
\newfont{\numerikNine}{cmss9}
\newfont{\numerikEight}{cmss8}
\newfont{\numerikSeven}{cmss7}
\newfont{\numerikSix}{cmss6}
\definecolor{Ema}{RGB}{0, 140, 250}
\begin{document}
\mainmatter              % start of a contribution
\title{Adaptive Time-Step Semi-Implicit One-Step Taylor Scheme for Stiff Ordinary Differential Equations}
%
%\titlerunning{Hamiltonian Mechanics}  % abbreviated title (for running head)
\author{S. Boscarino \and E. Macca*}
%
%\authorrunning{Ivar Ekeland et al.} % abbreviated author list (for running head)
%Craig Chambers, Kim B. Bruce, and Elisa Bertino}
%
\institute{University of Catania, Department of Mathematics and Computer Science, Catania 95127, Italy,\\
\email{sebastiano.boscarino@unict.it, emanuele.macca@unict.it},}

\maketitle              % typeset the title of the contribution

\begin{abstract}
In this study, we propose high-order implicit and semi-implicit schemes for solving ordinary differential equations (ODEs) based on Taylor series expansion. These methods are designed to handle stiff and non-stiff components within a unified framework, ensuring stability and accuracy. The schemes are derived and analyzed for their consistency and stability properties, showcasing their effectiveness in practical computational scenarios.
\keywords{Semi-implicit, high-order, Taylor expansion, ODE, time-step control}
\end{abstract}

\section{Introduction}
The numerical integration of ordinary differential equations (ODEs) is a fundamental task in scientific computing, with applications spanning physics, engineering, biology, and finance. The challenge often lies in dealing with stiff systems where explicit methods require prohibitively small time steps to maintain stability \cite{Wanner1996,Carrillo2021,Carrillo2023}. In contrast, implicit methods, while stable, can be computationally intensive due to the necessity of solving nonlinear systems at each time step \cite{Hundsdorfer2003}. 

To address these issues, semi-implicit (SI) and implicit-explicit (IMEX) techniques have gained prominence. These methods treat the non-stiff components explicitly \cite{Baeza2018,Baeza2020,Exner2024} and the stiff components implicitly, thus combining the stability advantages of implicit methods with the computational efficiency of explicit ones \cite{Boscarino2006,Scott2000}. High-order SI and IMEX schemes, in particular, offer enhanced accuracy and are valuable for long-time integration of stiff ODEs \cite{Boscarino2017,Baeza2017}.

There exist various approaches for implicit schemes based on Taylor series expansion. In this work, we develop high-order implicit and semi-implicit schemes based on Taylor expansion \cite{Kirlinger1991}. By employing the Taylor expansion, we can systematically construct schemes that are consistent to higher orders \cite{Miletics2004,Miletics2005}. We specifically focus on first and second-order schemes and analyze their stability properties through theoretical and numerical approaches. 

Consider a system of differential equations given by:
\begin{equation*}
U' = f(U) + g(U),
\end{equation*}
where $ U = U(t) \in \mathbb{R}^d $ is the vector of state variables, $f(U)$ represents the non-stiff part, and $g(U)$ represents the stiff part. 

However, explicit schemes are only feasible when the stiff part $g$ is not dominant. For strongly stiff problems, we develop semi-implicit schemes where the stiff part is treated implicitly to ensure stability.

In the following sections, we derive first and second-order semi-implicit Taylor schemes (SI-T-1, SI-T-2), analyze their stability properties, and demonstrate their application through numerical experiments. Our results indicate that these schemes are not only consistent and stable but also computationally efficient.

Moreover, the use of step-control techniques  \cite{Ahmed2015,Yao2014} and comparison with IMEX Runge-Kutta schemes \cite{Balac2013,Izzo2017} provide further insight into the quality of the solutions obtained. These techniques allow for adaptive time-stepping, ensuring both accuracy and stability, making our proposed methods robust and reliable for practical applications. Other time-step control methods, such as the \textit{a posteriori} Multi-dimensional Optimal Order Detection (MOOD) paradigm \cite{Loubère2024,Macca2024}, can be explored to assess the computational effort associated with time-step control. However, this is beyond the scope of the present work.

The rest of this paper is organized as follows. 
The second section introduces the governing equations, the derivation of the first and second order semi-implicit and implicit schemes. The third section presents the stability analysis of the schemes. The fourth section is related to the adaptive time-step control. A second order embedded IMEX Runge-Kutta scheme, for comparison with these new schemes, is also presented. Finally, numerical results are gathered in the fifth section to assess the good behavior of the semi-implicit methods. Conclusions are finally drawn.

\section{Numerical Schemes.}
Let us consider the system of differential equation
\begin{equation}
\label{gov_eq}
    U' = f(U) + g(U)
\end{equation} 
where $U=U(t):\mathbb{R}^{+}\rightarrow\mathbb{R}^d$ denotes the vector variable, $f(U)=f(U(t))$ the non-linear non-stiff term and $ g(U) = g(U(t))$ the non-linear stiff term.
We denote $U(0)=U_0$ the initial condition (IC) of the Cauchy problem
\begin{equation}
    \begin{cases}
        U' = f(U+g(U) \\
        U_0 = U(0).
    \end{cases}
\end{equation}

In the high order mining, the second derivative of $U$ can be expressed as
$$ U'' = J_f U' + J_g U' = (J_f + J_g)(f + g), $$
where, $J_f$ and $J_g$ represent, respectively $\partial f/ \partial U$ and $\partial g/\partial U.$

Using the Taylor series expansion, $U(t + k)$ can be expressed by:
\begin{equation}
    \label{Taylor_exp}
    U(t + k) = U(t) + kU'(t) + \frac{k^2}{2}U''(t) + \mathcal{P}(t,k)
\end{equation}
where $\mathcal{P}(t,k)$ is the remainder term.

The simplest approach to construct a numerical one-step scheme is to treat all components explicitly. The first and second-order schemes can be written as follows:
\begin{align}
    \label{Expl_or_1}
    U^{n+1} &= U^n + \Delta t (f^n + g^n); \\
    \label{Expl_or_2}
    U^{n+1} &= U^n + \Delta t (f^n + g^n) + \frac{\Delta t^2}{2} \Bigl( J_f^n + J_g^n \Bigr) \Bigl( f^n + g^n \Bigr).    
\end{align}
where $U^n \approx U(t^n) = U(t+n\Delta t)$, $\dt$ the time step, while $f^n =f(U^n)$, $g^n =g(U^n)$, $J_f^n$   and $J_g^n$ represent the gradient for $f$ and $g$ at time $t^n.$

Nevertheless, these explicit schemes can be useful when the stiff part $g$ is not too dominant or when the time step $\dt$ is chosen sufficiently small to maintain numerical stability. However, for strongly stiff problems, implicit or semi-implicit schemes are necessary to effectively handle numerical stability.

In this regard, efficient semi-implicit, imex and fully-implicit scheme are now considered.

\subsection{Semi-implicit one-step Taylor}
Starting from the second order Taylor expansion \eqref{Taylor_exp} the semi-implicit first and second order schemes can be derived. Indeed, treating $f,$ $J_f$ and $J_g$ explicitly and $g$ implicitly, the first order semi-implicit Taylor scheme (SI-T-1) becomes
\begin{equation}
    \label{SI-T_or_1}
    U^{n+1} = U^n + \Delta t (f^n + g^{n+1})
\end{equation}
where $g^{n+1} = g(U^{n+1})$. 

The second order semi-implicit Taylor scheme (SI-T-2) is written as 
\begin{equation}
    \label{SI-T_or_2}
    U^{n+1} = U^n + \Delta t (f^n + g^{n+1})  + \frac{\Delta t^2}{2} \Bigl( J_f^n + J_g^n \Bigr) \Bigl( f^n - g^{n+1} \Bigr). 
\end{equation}

These semi-implicit Taylor scheme are consistent up to first and second order, respectively. Indeed, if $g^{n+1}$ is replaced with his Taylor expansion of zero and first order, after some algebraic manipulation the explicit schemes \eqref{Expl_or_1}-\eqref{Expl_or_2} are recovered.

\subsection{Implicit one-step Taylor}
Following the same idea of the semi-implicit approach, the first and second order scheme can be derived treating all the terms implicitly. In particular, the first order method appears
\begin{equation}
    \label{I-T_or_1}
    U^{n+1} = U^n + \Delta t (f^{n+1} + g^{n+1});
\end{equation}
meanwhile, the second order scheme:
\begin{equation}
    \label{I-T_or_2}
    U^{n+1} = U^n + \Delta t (f^{n+1} + g^{n+1}) - \frac{\Delta t^2}{2} \Bigl( J_f^{n+1} + J_g^{n+1} \Bigr) \Bigl( f^{n+1} + g^{n+1} \Bigr)
\end{equation}
where $J_f^{n+1} $ and $ J_g^{n+1}$ represent the gradient for $f$ and $g$ at time $t^{n+1}.$

\section{Stability analysis}
In this section, the goal is to develop a theoretical stability analysis, following the standard approach used for linear ordinary differential equation (ODE) systems, i.e., $u' = \lambda u + \nu u$ with $\lambda,\,\nu \in \mathbb{C}$,  ${\rm Re}(\lambda,\nu)<0$ and $|\lambda|\ll |\nu|$. To simplify the discussion, we will focus on studying a scalar linear differential equation with $f$ an $g$ linear term \cite{CAT-SI,Hundsdorfer2003}.

Now the behaviour of the numerical solution for each first and second order semi-implicit, \eqref{SI-T_or_1} and \eqref{SI-T_or_2}, and implicit,\eqref{I-T_or_1} and \eqref{I-T_or_2}, scheme have been analyzed. In the limit case, ${\rm Re}(\nu) \to -\infty$, the exact solution is $u=0$, then we will show which scheme achieves correctly the exact solution. 

The SI-T-1 scheme applied to linear differential equation becomes:
\begin{equation*}
    u^{n+1} = u^{n} + \dt\Bigl(\lambda u^{n} + \nu u^{n+1}\Bigr) 
\end{equation*}
and defining $z = \dt\lambda$ and $w = \dt \nu,$ we obtain:
\begin{equation*}
    (1-w) u^{n+1} = (1+z)u^{n}.
\end{equation*}
Hence, we get the following numerical solution 
\begin{equation*}
    \label{equ:a_stable_func}
    u^{n+1} =\mathcal{R}_{\rm SI-T-1}(z,w) u^{n}
\end{equation*}
with 
\begin{align}
    \mathcal{R}_{\rm SI-T-1}(z,w) &= \frac{1+z}{1-w}.
\end{align}
In similar way 
\begin{align}
    \mathcal{R}_{\rm SI-T-2}(z,w) &= \frac{1+z + \ha(z^2 + zw) }{1-w + \ha(zw + w^2)};
    \\
    \mathcal{R}_{\rm I-T-1}(z,w) &= \frac{1}{1- z -w};
    \\
    \mathcal{R}_{\rm I-T-2}(z,w) &= \frac{1}{1-(z+w)(1-\ha(z + w))}.  
\end{align}
When $z\in\mathbb{C}^-$ and ${\rm Re}(w)\rightarrow-\infty$, $\mathcal{R}_{\rm all}(z,w) \rightarrow0$ and this implies that all the schemes based on the Taylor expansion proposed are $L-$stable.
\begin{figure}[!ht]
    \centering
     \includegraphics[width=0.85\textwidth]{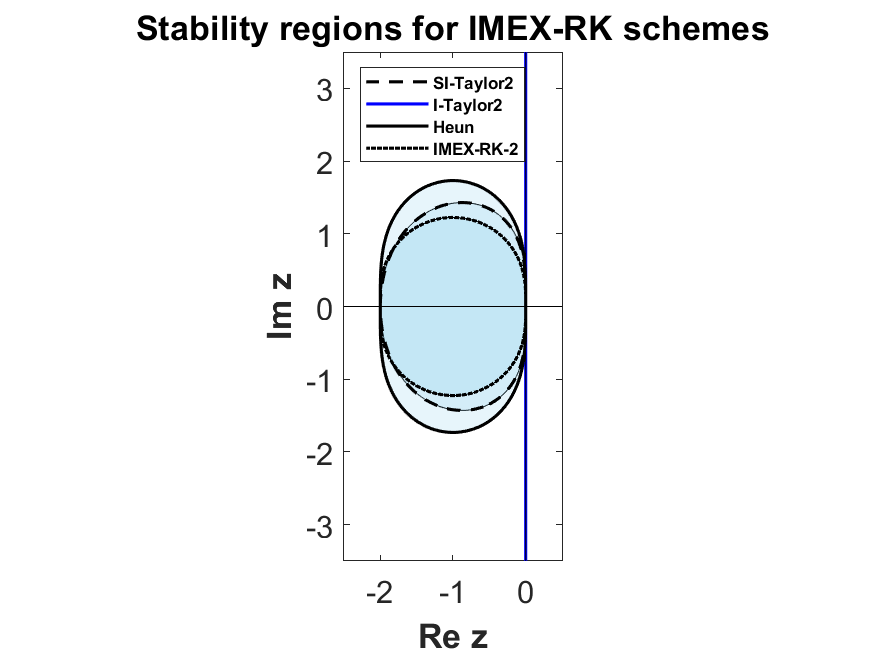}
     \caption{Stability regions $S_1$, in the $(Re(w); Im(w))$  plane, of some second order schemes such as: IMEX-RK2 method; Implicit and Semi-Implicit Taylor method and Heun scheme.}
     \label{Stab_or_2}
\end{figure}

Based on the analytical analysis of the linear stability of the schemes presented, it is evident that all the schemes are consistent. However, to facilitate a graphical comparison of different schemes, whether well-known in the literature or otherwise, we must first define the stability region for the given type of equation. This definition will enable us to compare the various schemes effectively.

Typically, the region of absolute stability $S$ associated with SI or IMEX-RK schemes is defined as
$$
S = \{(z, w) \in \mathbb{C}^2 : |\mathcal{R}(z,w)| \le 1\}.
$$
However, we gain a deeper understanding of the geometrical structure of the region $S$ associated with an SI or IMEX-RK scheme by demonstrating the existence of two regions in the complex plane, $S_1 \subset \mathbb{C}$ and $S_2 \subset \mathbb{C}$, with the following properties:
$$
S_1\times S_2 \subset S, \, \, \mathbb{C}^- \equiv \{w \in \mathbb{C}: \text{Re}(w) \le 0\} \subset S_2.
$$
These two sets, if they exist, are not unique. We consider two common choices to uniquely define the sets. Here, we present one choice by defining $S_1$ as the largest stability region for $z$ such that the scheme remains $A$-stable, i.e., $\mathbb{C}^- \subset S_2$.
This region is formally defined as
$$
S_1 = \{z \in \mathbb{C} : \sup_{w \in \mathbb{C}} |\mathcal{R}(z,w)| \le 1\}.
$$
{For more details on the other choice, the definition of $S_2$ and the computation of $S_1$, refer to the book \cite{IMEX_book}. In Figure~\ref{Stab_or_2}, we show the stability regions $S_1$ for some SI and IMEX-RK schemes introduced in the previous section. It depicts the stability region $S_1$ for the second-order schemes, including: IMEX-RK2 \cite{Boscarino2016,BUMI2024,Exner2024}, Implicit and Semi-Implicit Taylor2 (I-T2 and SI-T2), and Heun. } 

{Since the two semi-implicit and implicit schemes are designed such that $\mathcal{R}(z, w) \rightarrow 0$ as $\text{Re}(w) \rightarrow -\infty$, they are $L-$stable \cite{Wanner1996}. 
The $L-$stability is strongly related to the asymptotic preserving (AP) property, especially if they are applied in the context of PDEs, \cite{IMEX_book}.}

{Furthermore, it is straightforward to show that in the limit case $\text{Re}(w) \rightarrow -\infty$, the schemes perform well for both well-prepared initial conditions and not well-prepared ones\footnote{A well-prepared initial condition in the context of a numerical simulation or mathematical model refers to an initial state that is carefully selected to reflect the problem being studied and to ensure the stability and accuracy of the simulation. It means that the initial condition is chosen in a way that is physically meaningful and conducive to obtaining reliable results from the simulation. For instance, in the Van der Pol model, given $y(0)$, the initial condition is well-prepared if and only if $z(0) = y(0)/(1-y(0)^2)$.}.}

\section{Time-step controller}
Sometimes, it is useful to vary the time step size of a numerical method to solve difficult problems, such as stiff problems. To achieve this, it is important to choose a time-step controller that ensures both accuracy and stability. This can be done by estimating and controlling some measure of the local error, as extensively explained in \cite{Wanner1996}. The fundamental idea behind a time-step controller is usually to define an embedded time-integration method based on the main numerical method of order $p$. This involves providing an additional scheme, called the embedded method.  %that is either one order lower $(q=p-1)$ or one order higher $(q=p+1)$ than the main scheme to compute the local error.
An embedded RK method associated to a RK one of order $p$, is a scheme with the same matrix $A$, and nodes $c$, and with a new set of weights $\hat{b}$, computed imposing that the embedded scheme has order $q = p-1$ or $q = p + 1$, i.e., one order less (or more) accurate than the main scheme. These embedded methods are designed to produce an estimate of the local error for a single Runge–Kutta step and are used to control the local error for the adaptive time step controller.

%\ema{Quali di questi due stiamo usando? Il primo?}

% "(2,1)-IMEX-RK" scheme is a second order IMEX-RK coupled with a first order IMEX-RK one, i.e.
% \[
% \begin{array}{c|ccc}
%  0 & 0 & 0 & 0\\
% \frac{2}{5} & \frac{2}{5} & 0 & 0\\
%  1 & 0 & 1 & 0 \\
% \hline
%  &0& \frac{5}{6} & \frac{1}{6}\\[+0.1cm]
%  \hline
%  &0& \frac{4}{5} & \frac{1}{5} \\
% \end{array}
% \quad \quad 
% \begin{array}{c|ccc}
%  0 & 0 & 0 & 0\\
% \frac{2}{5} & 0 & \frac{2}{5} & 0 \\
% 1 & 0 & \frac{5}{6} & \frac{1}{6} \\
% \hline
% &0 & \frac{5}{6} & \frac{1}{6}\\[+0.1cm]
% \hline
% & 0 & \frac{4}{5} & \frac{1}{5}
% \end{array}
% \]
An example is the (2,1)-DIRK scheme, which is a second-order implicit DIRK scheme coupled with a first-order implicit one:
\begin{equation}\label{A1}
\begin{array}{c|cc}
 \gamma & \gamma & 0\\
 1 & 1-\gamma  & \gamma \\
\hline
  & 1-\gamma  & \gamma \\[+0.1cm]
  \hline
& 0 & 1
\end{array} \, \, .
\end{equation}
In the numerical tests, we combine this (2,1)-DIRK scheme with an explicit one:
\begin{equation}\label{A2}
\begin{array}{c|cc}
 \gamma & 0 & 0\\
 1 & c  & 0 \\
\hline
  & 1-\gamma  & \gamma \\[+0.1cm]
  \hline
& 0 & 1
\end{array}
\end{equation}
with $\gamma = 1-\sqrt{2}/2$ and $c = 1/(2\gamma)$. We call this IMEX-RK(2,1) scheme.

Computing the local error can be useful for automatically and adaptively controlling the time step $\Delta t$. Several well-known controllers in the literature, such as the I, PI, or PID controllers, facilitate this process (see \cite{Wanner1996} for details). The two most common controllers are the I and PID controllers. Here we consider the I controller defined by 
\begin{equation}\label{Icontr}
(\Delta t)_I^{n+1} = \kappa \Delta t^n \left(\frac{Tol}{\| \delta^n\|_{\infty}} \right)^{1/q}
\end{equation}
with local error estimate $\delta^n = U^n - \hat{U}^n$, where $\hat{U}^n$ and $q$ are the numerical solution
and the order of accuracy of the embedded method with $\kappa$ a safety factor to ensure success on the next try. Usually $\kappa$ is choose between $0.89$, $0.9$.; $Tol$ is a user specified
tolerance, $\Delta t$ is the step size of the last completed step.

%\subsection{Embedded second order Runge-Kutta}
%An embedded RK method associated to a RK one of order $p$, is a scheme with the same matrix $A$, and nodes $c$, and with a new set of weights $\hat{b}$, computed imposing that the embedded scheme has order $q = p-1$ or $q = p + 1$, i.e., one order less (or more) accurate than the main scheme. They are designed to produce an estimate of the local error for a single Runge–Kutta step, and are used to control the local error for the adaptive time step controller.

\section{Numerical experiment}
In this example, we analyze the behavior of several schemes: SI-T-1, SI-T-2 and IMEX-RK(2,1), when applied to the Van der Pol’s (VdP) problem
\begin{equation}
\label{VDP}
    \begin{cases}
        y' = z, \\ z' = \mu(1-y^2)z-y,        
    \end{cases}
\end{equation}
with well-prepared initial conditions 
\begin{equation}\label{IC_wp}
    {\rm IC_1} = \begin{cases}
        y(0) = 2 \\ z(0) = \frac{y(0)}{1 - y(0)^2}
    \end{cases}
\end{equation} 
so that no initial layer appears at the beginning, and unprepared initial conditions
\begin{equation}\label{IC_nowp}
    {\rm IC_2} = \begin{cases}
        y(0) = 2 \\ z(0) = 0.
    \end{cases}
\end{equation}
In this case, $$ U = \begin{bmatrix}
    y \\ z
\end{bmatrix} \quad\quad f(U) = \begin{bmatrix}
    z\\-y
\end{bmatrix} \quad \quad g(U) = \mu\begin{bmatrix}
    0 \\ (1-y^2)z
\end{bmatrix}.   $$
% 
% WP
\begin{figure}[!ht]
     \centering
     \begin{subfigure}[b]{0.33\textwidth}
         \centering
         \includegraphics[width=\textwidth]{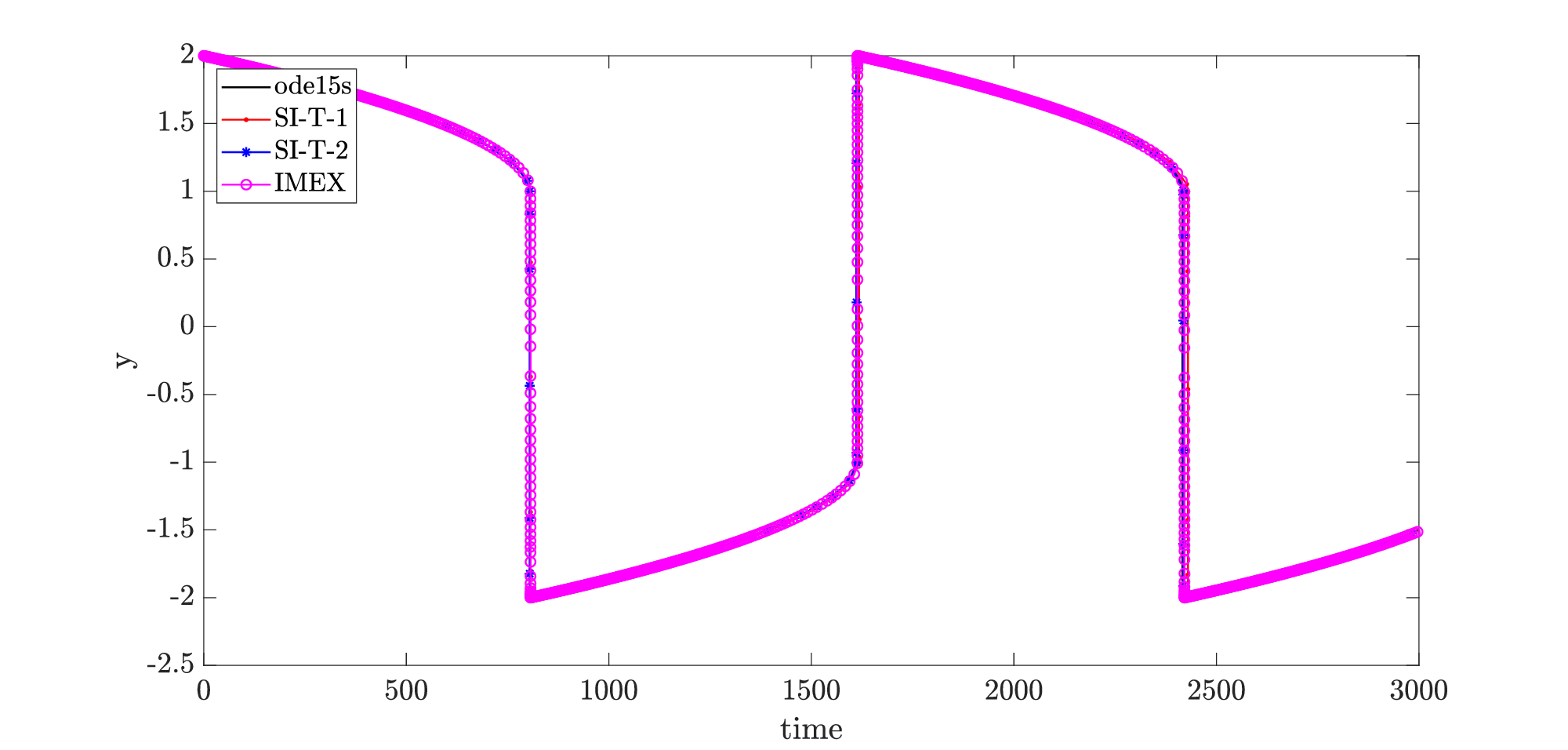}
         \caption{Numerical solution $y$.}
         \label{WP_y_SI}
     \end{subfigure}
     \hfill
     % \begin{subfigure}[b]{0.49\textwidth}
     %     \centering
     %     \includegraphics[width=\textwidth]{Plots/WP/SI-T-IMEX_z.eps}
     %     \caption{Numerical solution $z$.}
     %     \label{WP_z_SI}
     % \end{subfigure} 
     % \\
     \begin{subfigure}[b]{0.32\textwidth}
         \centering
         \includegraphics[width=\textwidth]{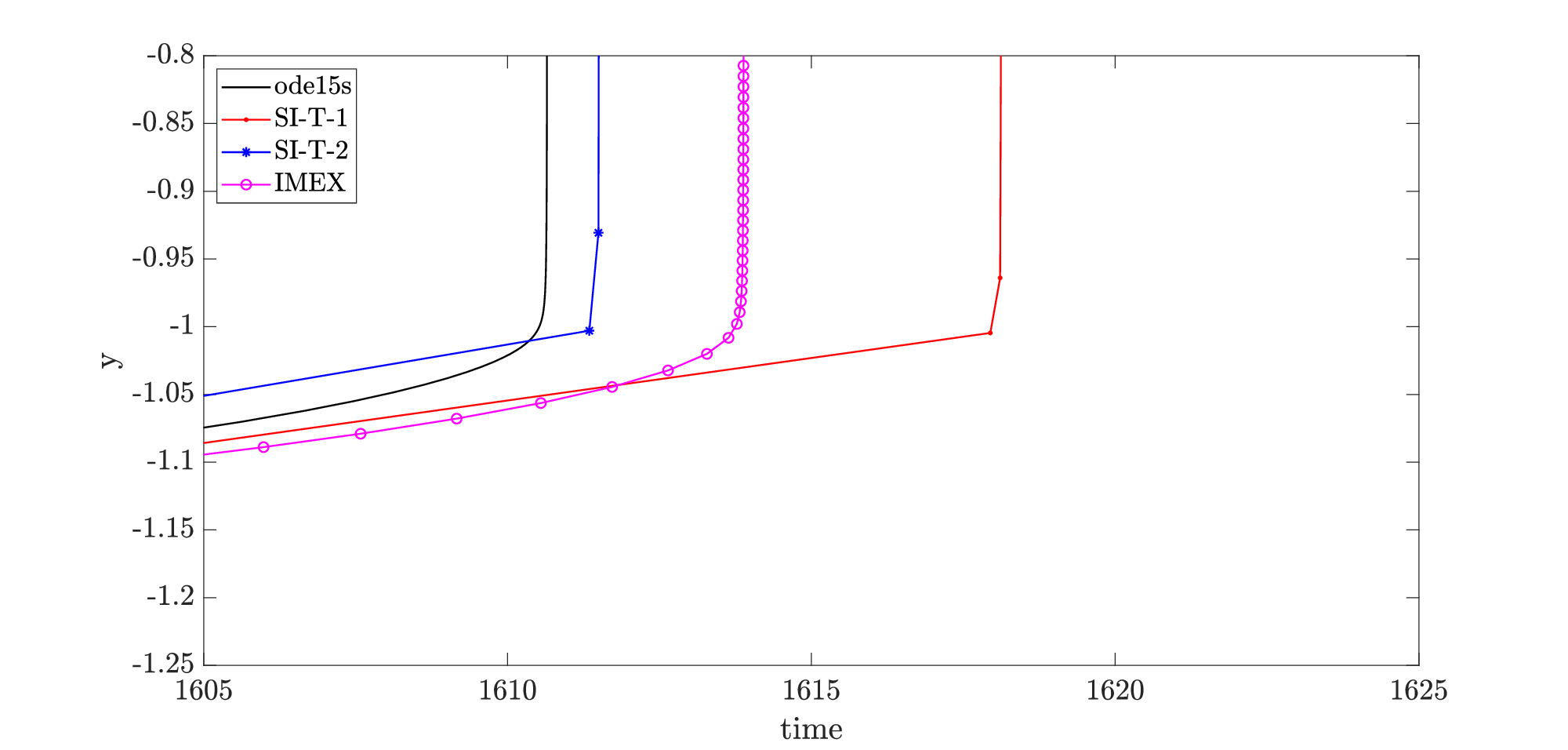}
         \caption{$y$ zoom bottom at $t \approx 1.6\mu$.}
         \label{WP_y_1_SI}
     \end{subfigure}
     \hfill
     % \begin{subfigure}[b]{0.49\textwidth}
     %     \centering
     %     \includegraphics[width=\textwidth]{Plots/WP/SI-T-IMEX_z_zoom_2.eps}
     %     \caption{$z$ zoom bottom at $t \approx 1.6\mu$.}
     %     \label{WP_z_1_SI}
     % \end{subfigure}
     \begin{subfigure}[b]{0.32\textwidth}
         \centering
         \includegraphics[width=\textwidth]{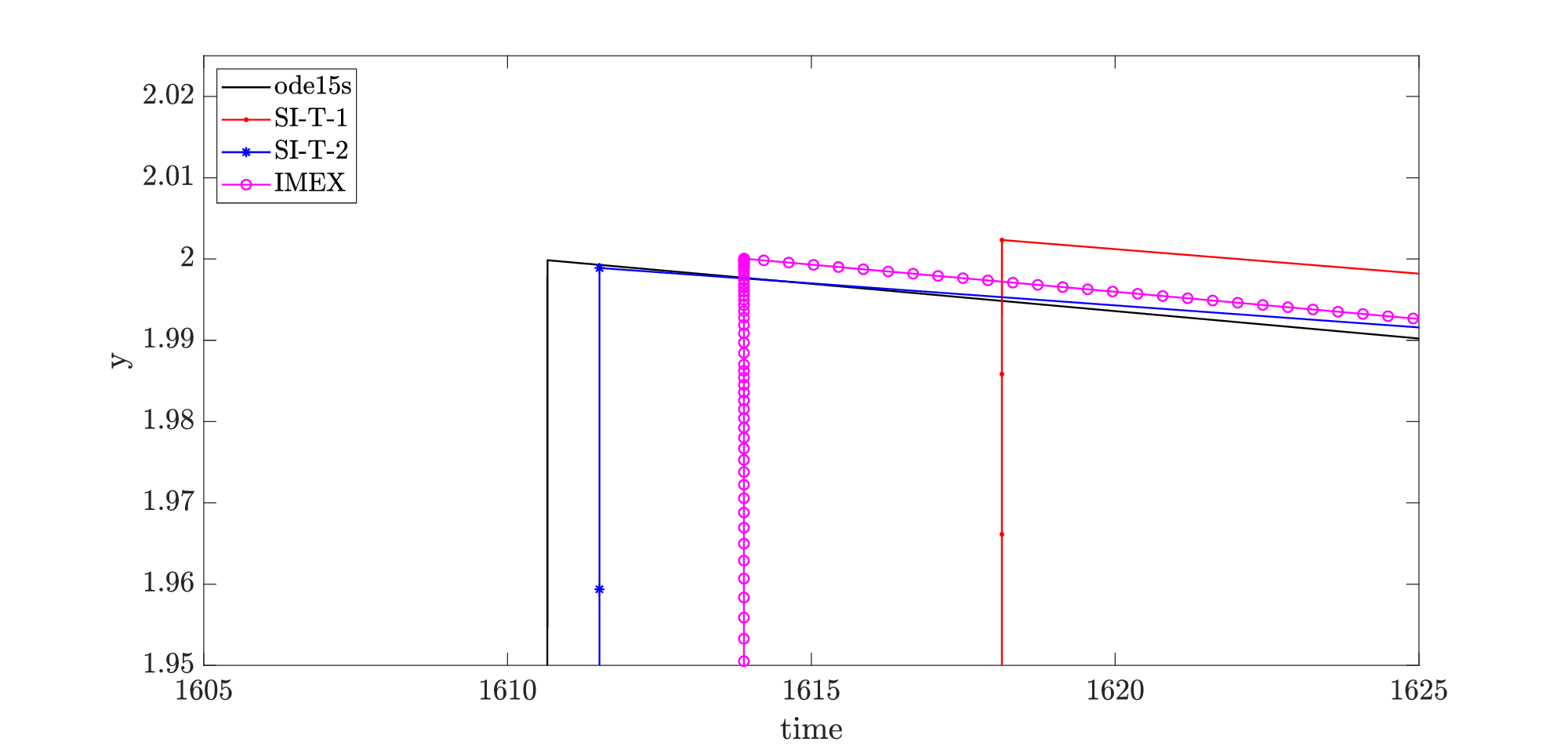}
         \caption{$y$ zoom top at $t \approx 1.6\mu$.}
         \label{WP_y_2_SI}
     \end{subfigure}
     % \hfill
     % \begin{subfigure}[b]{0.49\textwidth}
     %     \centering
     %     \includegraphics[width=\textwidth]{Plots/WP/SI-T-IMEX_z_zoom_1.eps}
     %     \caption{$z$ zoom top at $t \approx 1.6\mu$.}
     %     \label{WP_z_2_SI}
     % \end{subfigure}
     \caption{Semi-implicit numerical $y-$solution for the Van der Pol system \eqref{VDP} with well-prepared initial conditions \eqref{IC_wp} obtained at time $t=3\mu$ with $\dt_0 = 0.1$. The zoom of the second challenging boundary layers $t \approx 1.6\mu$ are reported. The reference solutions have been computed with the ode15s matlab solver. }
     \label{WP_SI}
\end{figure}
\begin{figure}[!ht]
     \centering
     \begin{subfigure}[b]{0.32\textwidth}
         \centering
         \includegraphics[width=\textwidth]{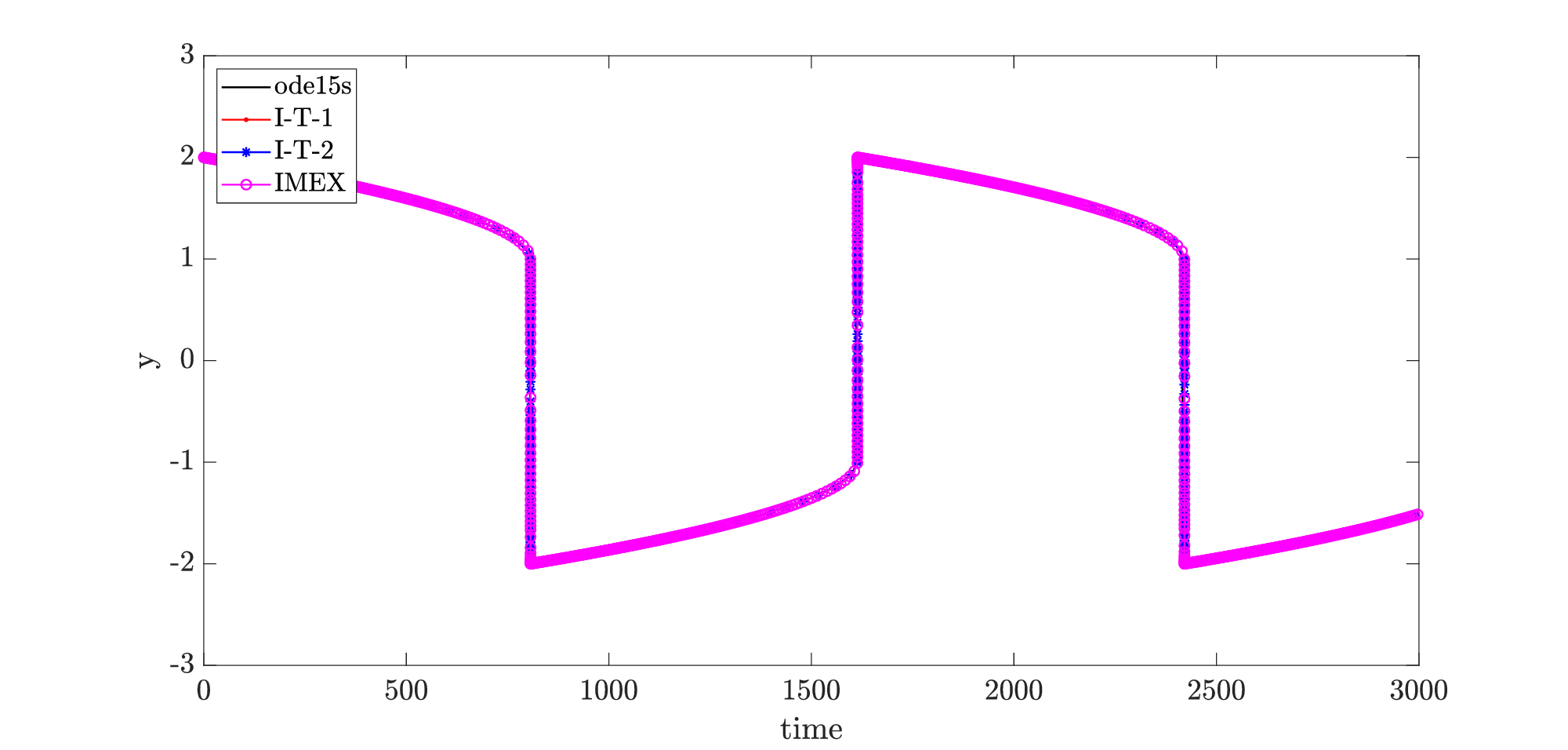}
         \caption{Numerical solution $y$.}
         \label{WP_y_I}
     \end{subfigure}
     \hfill
     % \begin{subfigure}[b]{0.49\textwidth}
     %     \centering
     %     \includegraphics[width=\textwidth]{Plots/WP/I-T-IMEX_z.eps}
     %     \caption{Numerical solution $z$.}
     %     \label{WP_z_I}
     % \end{subfigure} 
     % \\
     \begin{subfigure}[b]{0.32\textwidth}
         \centering
         \includegraphics[width=\textwidth]{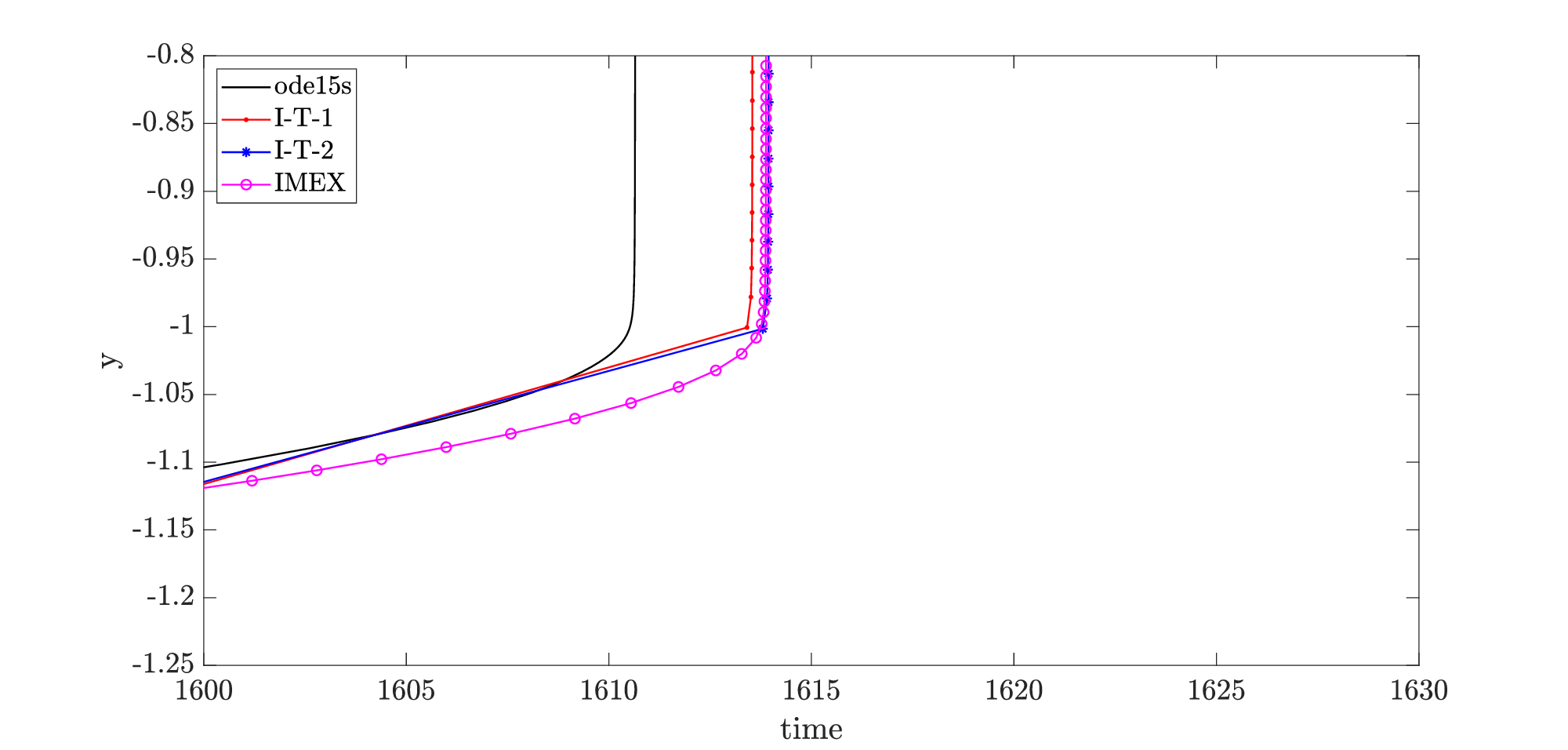}
         \caption{$y$ zoom bottom at $t \approx 1.6\mu$.}
         \label{WP_y_1_I}
     \end{subfigure}
     \hfill
     % \begin{subfigure}[b]{0.49\textwidth}
     %     \centering
     %     \includegraphics[width=\textwidth]{Plots/WP/I-T-IMEX_z_zoom_2.eps}
     %     \caption{$z$ zoom bottom at $t \approx 1.6\mu$.}
     %     \label{WP_z_1_I}
     % \end{subfigure}
     \begin{subfigure}[b]{0.32\textwidth}
         \centering
         \includegraphics[width=\textwidth]{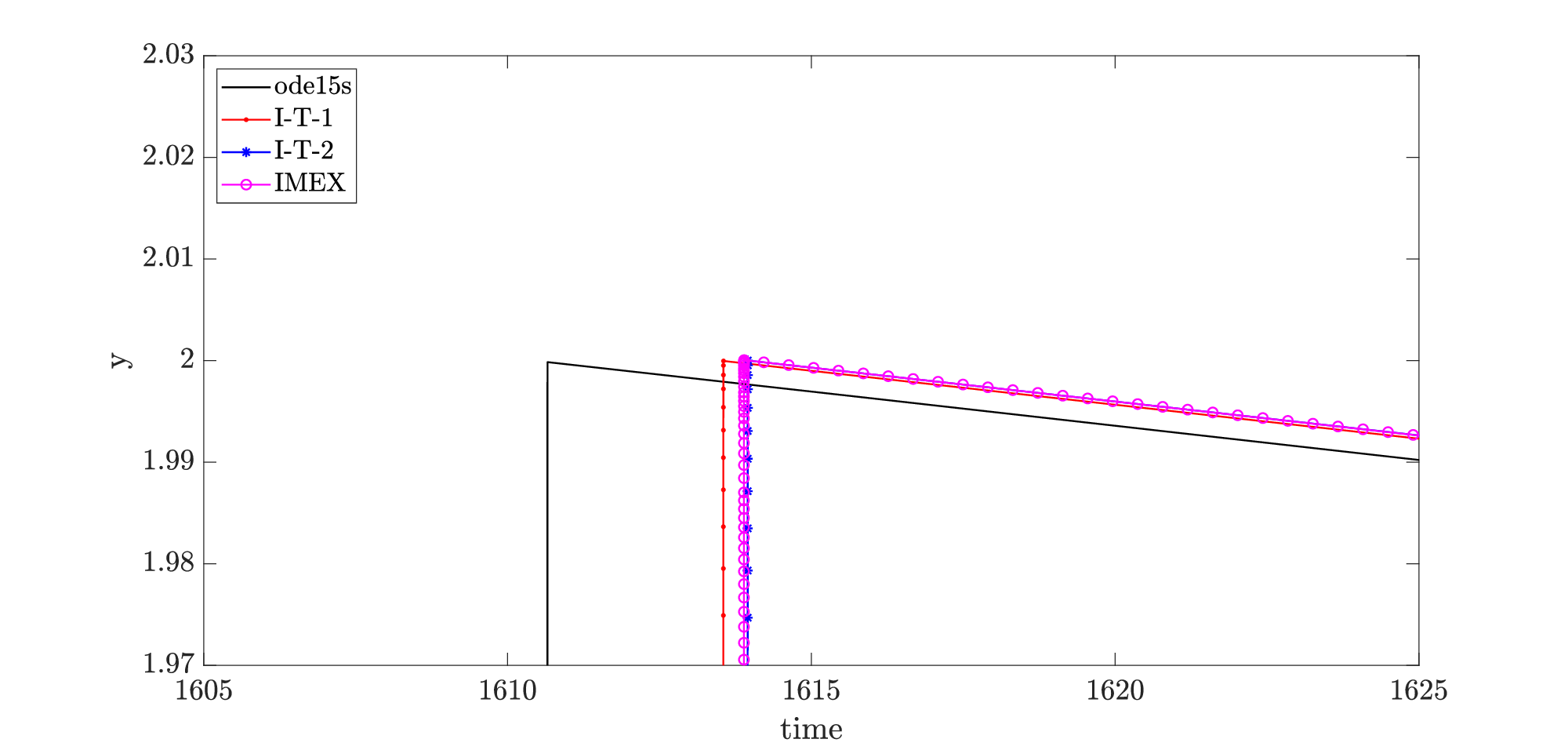}
         \caption{$y$ zoom top at $t \approx 1.6\mu$.}
         \label{WP_y_2_I}
     \end{subfigure}
     % \hfill\begin{subfigure}[b]{0.49\textwidth}
     %     \centering
     %     \includegraphics[width=\textwidth]{Plots/WP/I-T-IMEX_z_zoom_1.eps}
     %     \caption{$z$ zoom top at $t \approx 1.6\mu$.}
     %     \label{WP_z_2_I}
     % \end{subfigure}
     \caption{Implicit numerical $y-$solutions for the Van der Pol system \eqref{VDP} with well-prepared initial conditions \eqref{IC_wp} obtained at time $t=3\mu$ with $\dt_0 = 0.1$. The zoom of the second challenging boundary layers $t \approx 1.6\mu$ are reported. The reference solutions have been computed with the ode15s matlab solver. }
     \label{WP_I}
\end{figure}
% 
% noWP
\begin{figure}[!ht]
     \centering
     \begin{subfigure}[b]{0.32\textwidth}
         \centering
         \includegraphics[width=\textwidth]{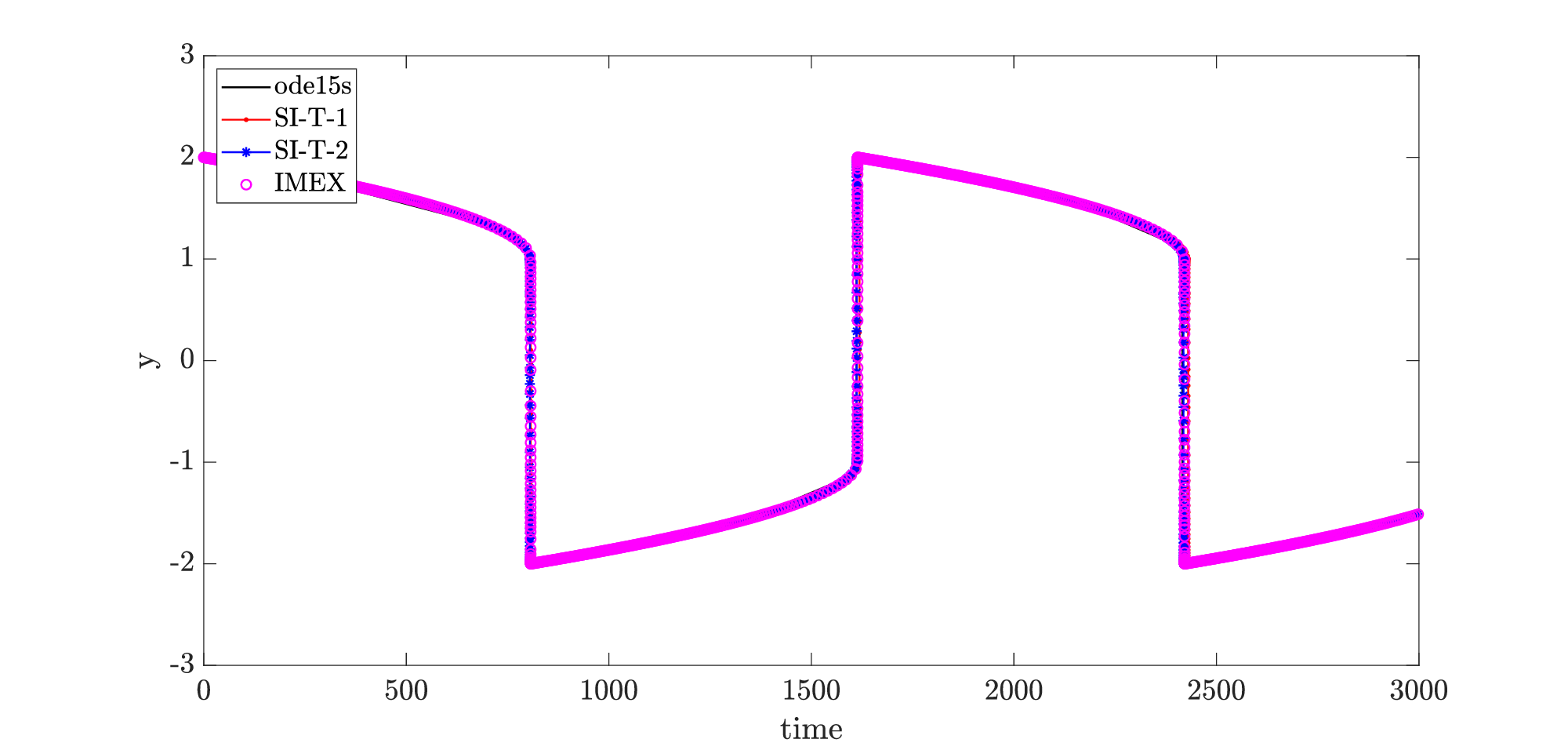}
         \caption{Numerical solution $y$.}
         \label{noWP_y_SI}
     \end{subfigure}
     \hfill
     % \begin{subfigure}[b]{0.49\textwidth}
     %     \centering
     %     \includegraphics[width=\textwidth]{Plots/noWP/SI-T-IMEX_z.eps}
     %     \caption{Numerical solution $z$.}
     %     \label{noWP_z_SI}
     % \end{subfigure} 
     % \\
     \begin{subfigure}[b]{0.32\textwidth}
         \centering
         \includegraphics[width=\textwidth]{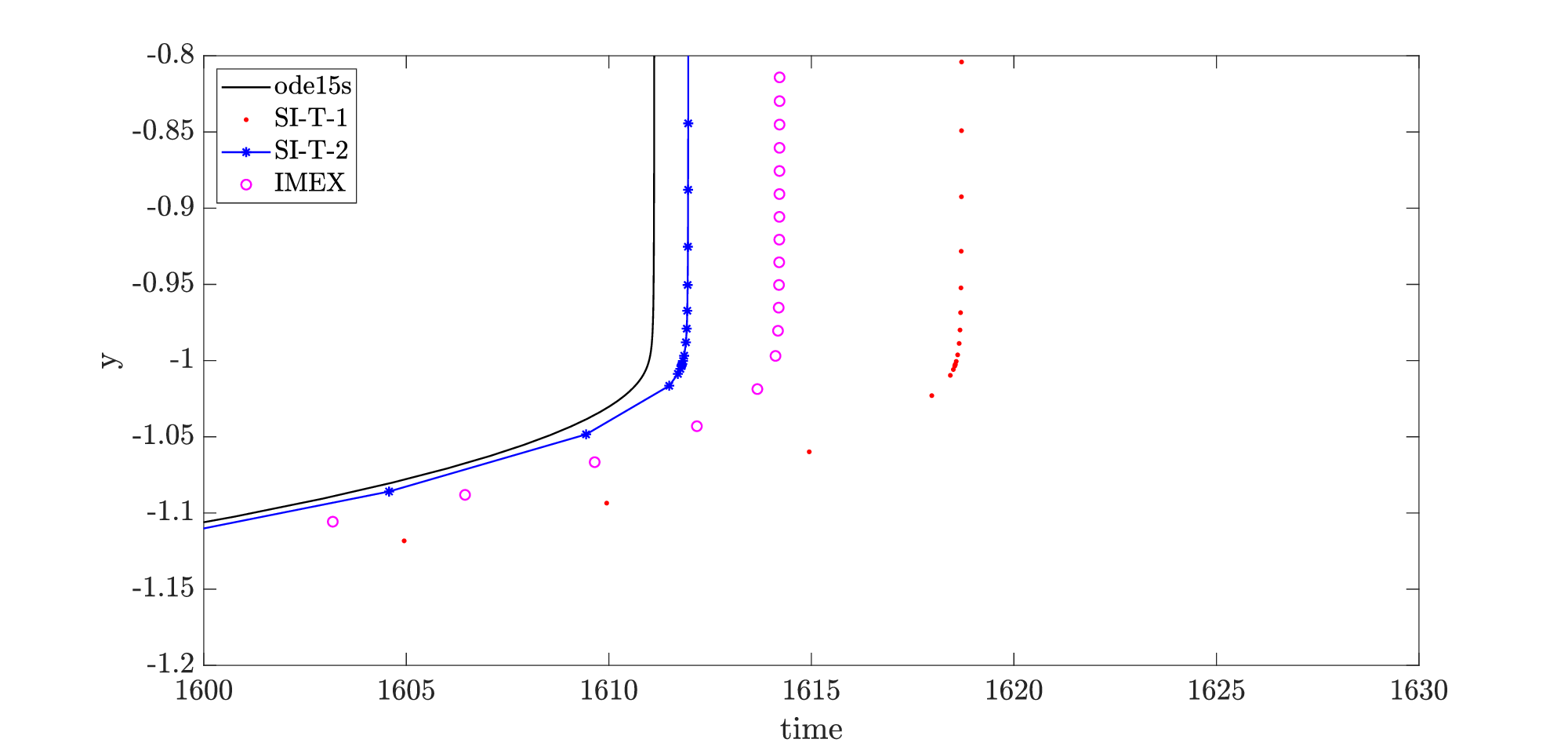}
         \caption{$y$ zoom bottom at $t \approx 1.6\mu$.}
         \label{noWP_y_1_SI}
     \end{subfigure}
     \hfill
     % \begin{subfigure}[b]{0.49\textwidth}
     %     \centering
     %     \includegraphics[width=\textwidth]{Plots/noWP/SI-T-IMEX_z_zoom_2.eps}
     %     \caption{$z$ zoom bottom at $t \approx 1.6\mu$.}
     %     \label{noWP_z_1_SI}
     % \end{subfigure}
     \begin{subfigure}[b]{0.32\textwidth}
         \centering
         \includegraphics[width=\textwidth]{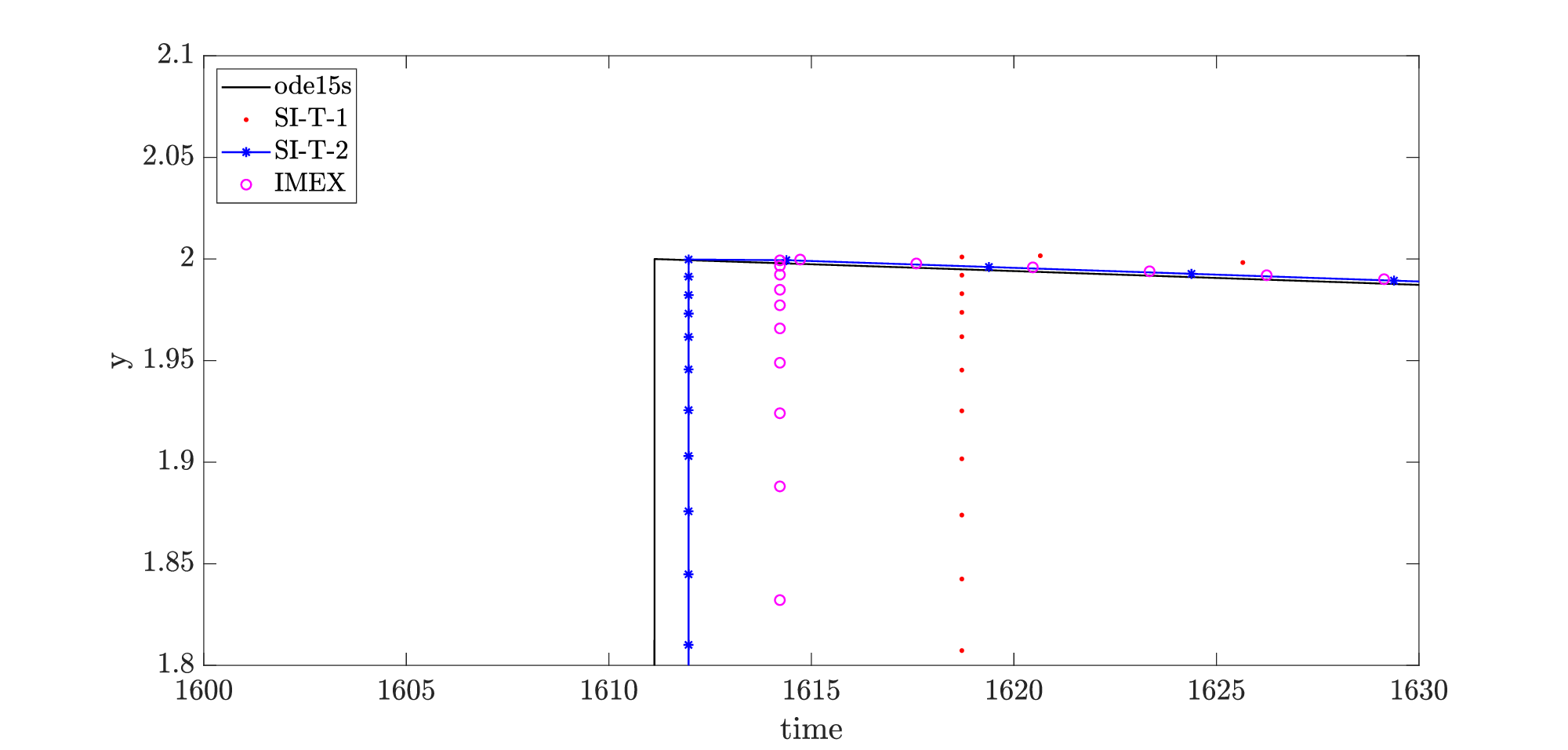}
         \caption{$y$ zoom top at $t \approx 1.6\mu$.}
         \label{noWP_y_2_SI}
     \end{subfigure}
     % \hfill\begin{subfigure}[b]{0.49\textwidth}
     %     \centering
     %     \includegraphics[width=\textwidth]{Plots/noWP/SI-T-IMEX_z_zoom_1.eps}
     %     \caption{$z$ zoom top at $t \approx 1.6\mu$.}
     %     \label{noWP_z_2_SI}
     % \end{subfigure}
     \caption{Semi-implicit numerical $y-$solution for the Van der Pol system \eqref{VDP} with unprepared initial conditions \eqref{IC_nowp} obtained at time $t=3\mu$ with $\dt_0 = 0.1$. The zoom of the second challenging boundary layers $t \approx 1.6\mu$ are reported. The reference solutions have been computed with the ode15s matlab solver. }
     \label{noWP_SI}
\end{figure}
\begin{figure}[!ht]
     \centering
     \begin{subfigure}[b]{0.32\textwidth}
         \centering
         \includegraphics[width=\textwidth]{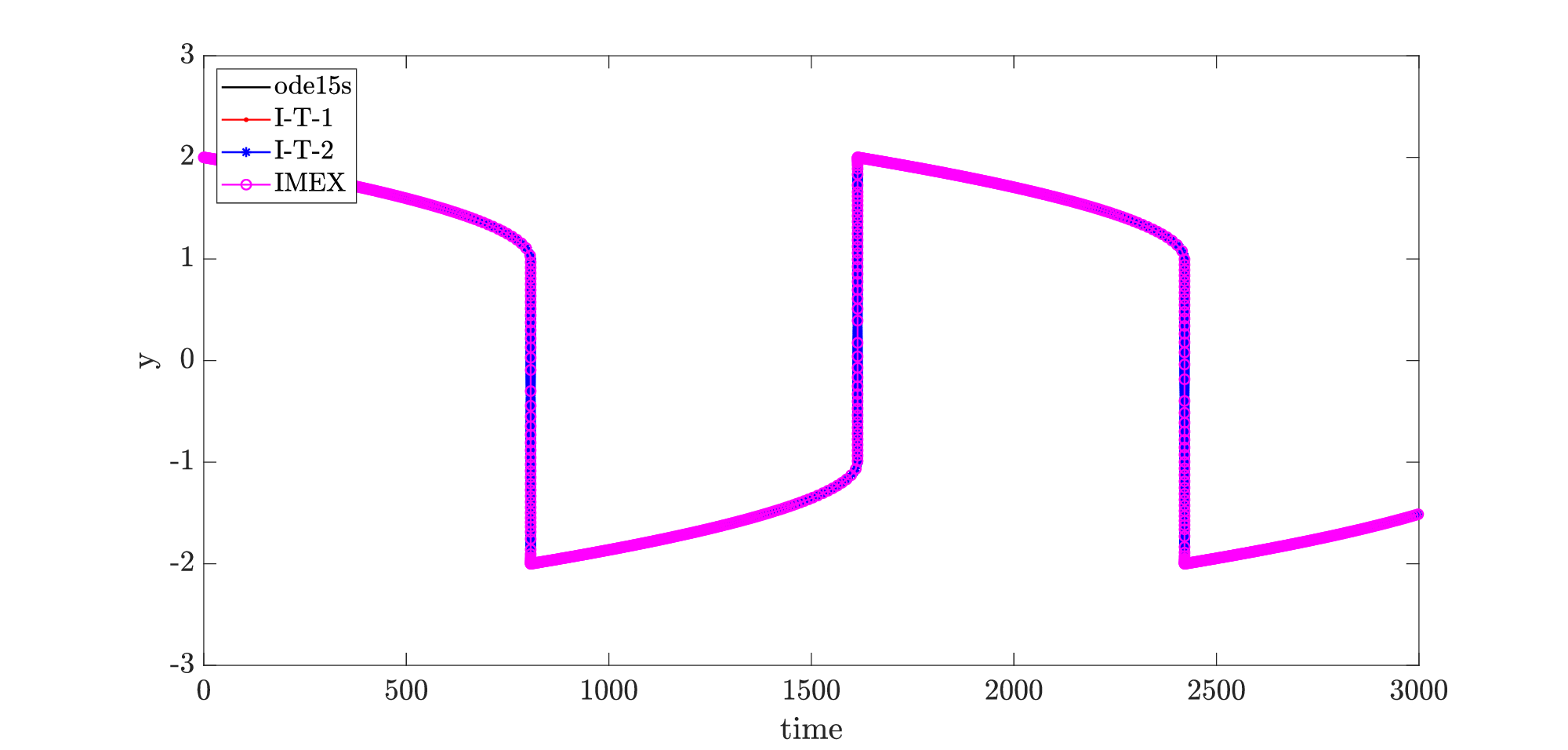}
         \caption{Numerical solution $y$.}
         \label{noWP_y_I}
     \end{subfigure}
     \hfill
     % \begin{subfigure}[b]{0.49\textwidth}
     %     \centering
     %     \includegraphics[width=\textwidth]{Plots/noWP/I-T-IMEX_z.eps}
     %     \caption{Numerical solution $z$.}
     %     \label{noWP_z_I}
     % \end{subfigure} 
     % \\
     \begin{subfigure}[b]{0.32\textwidth}
         \centering
         \includegraphics[width=\textwidth]{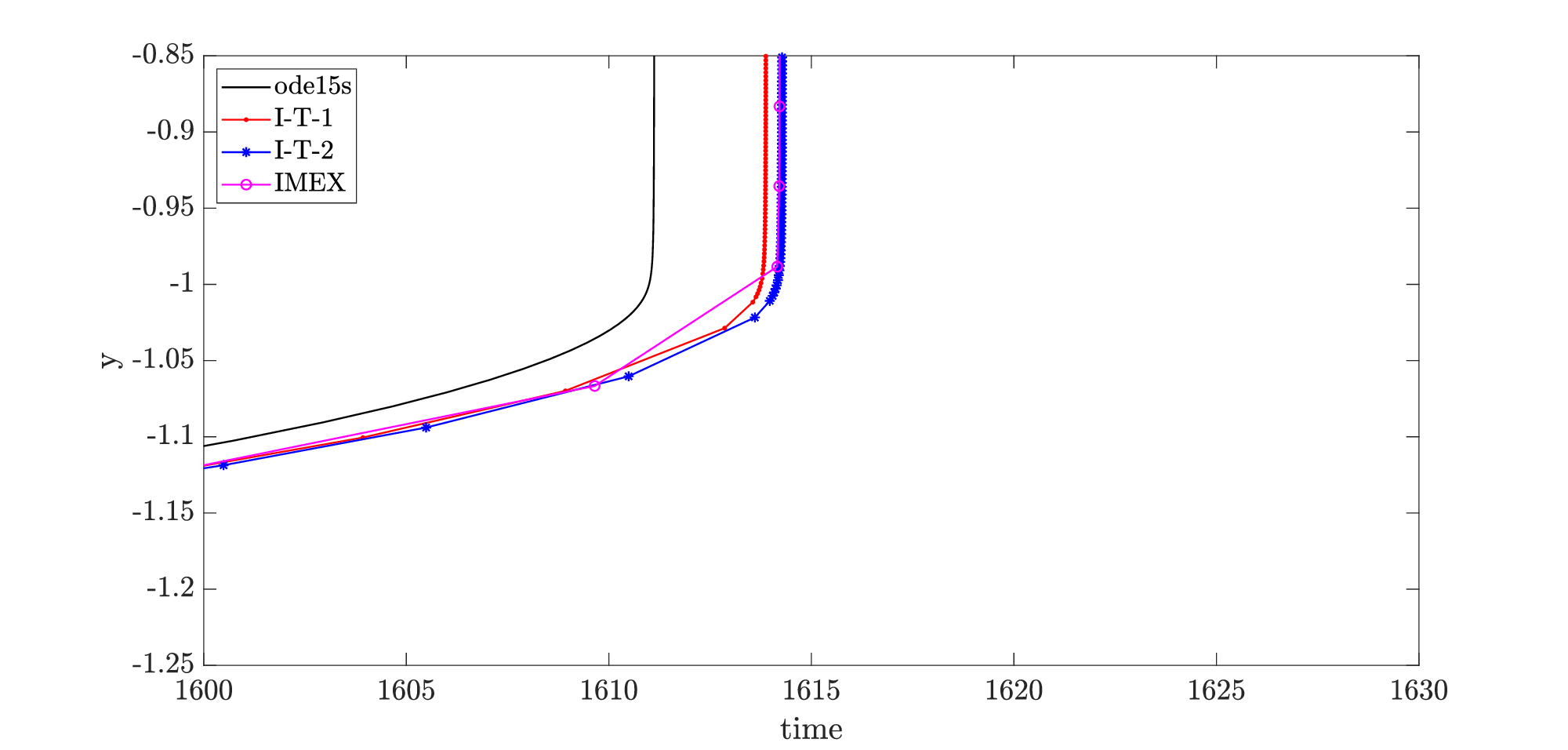}
         \caption{$y$ zoom bottom at $t \approx 1.6\mu$.}
         \label{noWP_y_1_I}
     \end{subfigure}
     \hfill
     % \begin{subfigure}[b]{0.49\textwidth}
     %     \centering
     %     \includegraphics[width=\textwidth]{Plots/noWP/I-T-IMEX_z_zoom_2.eps}
     %     \caption{$z$ zoom bottom at $t \approx 1.6\mu$.}
     %     \label{noWP_z_1_I}
     % \end{subfigure}
     \begin{subfigure}[b]{0.32\textwidth}
         \centering
         \includegraphics[width=\textwidth]{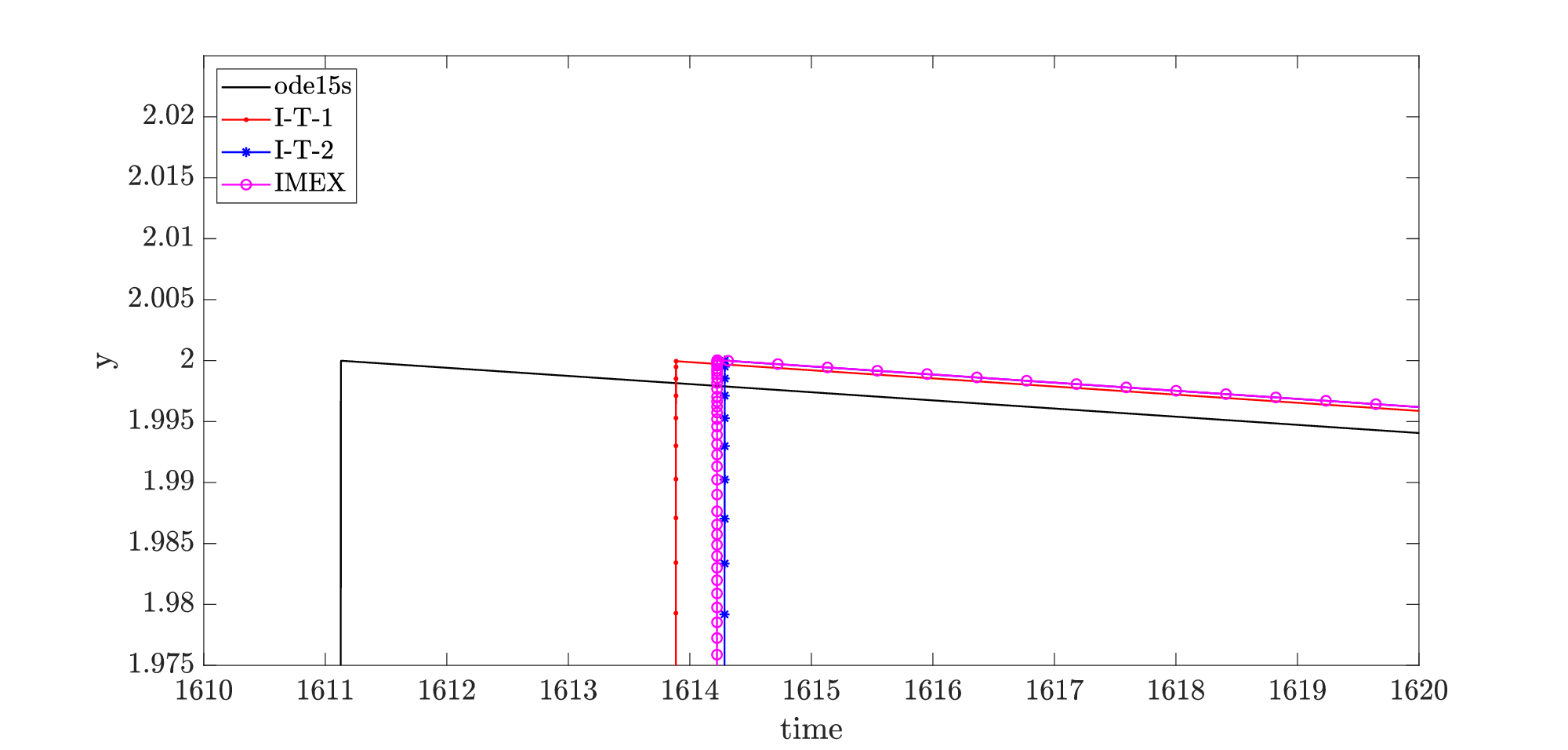}
         \caption{$y$ zoom top at $t \approx 1.6\mu$.}
         \label{noWP_y_2_I}
     \end{subfigure}
     % \hfill\begin{subfigure}[b]{0.49\textwidth}
     %     \centering
     %     \includegraphics[width=\textwidth]{Plots/noWP/I-T-IMEX_z_zoom_1.eps}
     %     \caption{$z$ zoom top at $t \approx 1.6\mu$.}
     %     \label{noWP_z_2_I}
     % \end{subfigure}
     \caption{Implicit numerical $y-$solution for the Van der Pol system \eqref{VDP} with unprepared initial conditions \eqref{IC_nowp} obtained at time $t=3\mu$ with $\dt_0 = 0.1$. The zoom of the second challenging boundary layers $t \approx 1.6\mu$ are reported. The reference solutions have been computed with the ode15s matlab solver. }
     \label{noWP_I}
\end{figure}
% --------------
% 
%
We set $\mu = 10^{3}$ and integrate over the interval $0 < t < 3\mu$. The VdP equation, when the parameter $\mu\gg1$, provides a challenging test because it develops very challenging boundary layers around the times $t \approx 0.8\mu$, $t \approx 1.6\mu$ and $t \approx 2.4\mu$. Consequently, a practical error controller (a reliable step-size controller) is required for the current method in order to compute the solution accurately. In our numerical test, we use the I-controller (\ref{Icontr}) to select the correct time step. We set $\kappa = 0.9$ as a safety factor. We applied these schemes with a tolerance $Tol =  10^{-5}$ and an initial step size of $\Delta t = 10^{-2}$, while $q=2$.

Figures~\ref{WP_SI}-\ref{WP_I} present the semi-implicit and implicit numerical solutions for the Van der Pol equations \eqref{VDP} with well-prepared initial conditions \eqref{IC_wp} at the final time $t = 3\mu$, using adaptive time-step control and an initial time-step $\dt_0 = 0.1$. The reference solutions were computed using the MATLAB ode15s solver, and the implicit and semi-implicit numerical solutions were compared with the embedded second-order IMEX RK(2,1) (\ref{A1})-(\ref{A2}) solutions\footnote{In the plots, to ensure the readability of the graphs, not all time-steps have been displayed. An appropriate selection, uniformly scaled, was considered, see Table~\ref{CPU_number_step}.}. Figures~\ref{noWP_SI}-\ref{noWP_I} display analogous results for the Van der Pol equations \eqref{VDP} with unprepared initial conditions \eqref{IC_nowp}. Specifically, detailed examinations of the solutions at time $t \approx 1.6\mu$ were conducted to highlight the differences between the various schemes. As observed, there are no significant differences between the results obtained from well-prepared and unprepared initial conditions. This behavior is due to the AP nature of the schemes considered.  

{Table \ref{CPU_number_step} presents the CPU time in seconds and the number of time steps required for both well-prepared and unprepared initial conditions across all schemes. The data clearly shows that the semi-implicit schemes are computationally more efficient, needing fewer iterations compared to other methods. For completeness, we also include the CPU time and the number of points for the ode15s solver. This MATLAB solver, based on the BDF linear multistep method, is computationally optimized for this type of equation. In contrast, the semi-implicit, implicit, and IMEX Runge-Kutta (2,1) schemes are not computationally optimized to the same extent, which impacts their performance.}
% ----- Tab ----------
\begin{table}[!ht]
%\numerikNine
\centering
\begin{tabular}{|c||c|c|c|c|c|c|}
\hline \multicolumn{7}{|c|}{\textbf{Van der Pol well-prepared}} \\
\hline  
\hline 
           &  ode15s & SI-T-1 & SI-T-2 & I-T-1 & I-T-2 & IMEX-RK(2,1) \\ \hline
CPU        & 0.05    & 10.76  & 10.75  & 26.88 & 26.98 & 16.65 \\
time-step  &  632    & 38602  & 38572  &160083 & 160103& 2819271  \\
\hline
\hline \multicolumn{7}{|c|}{\textbf{Van der Pol unprepared}} \\
\hline
\hline
CPU        & 0.05    & 10.96  &10.88   &25.96  & 25.99 & 16.48       \\
time-step  &  592    & 38547  &38563   &159692 & 159710& 2804550        \\
\hline
\end{tabular}
\vspace{0.2cm}
\caption{Van der Pol system \eqref{VDP} with well-prepared \eqref{IC_wp} and unprepared \eqref{IC_nowp}. CPU time in second and number of time-steps for both cases and all schemes.}
\label{CPU_number_step}
\end{table}

\section{Conclusion}
In this work, we developed and analyzed the semi-implicit schemes based on Taylor expansion for solving ordinary differential equations (ODEs) with stiff and non-stiff components. Specifically, we focused on first and second-order schemes and compared their performance against second order embedded IMEX Runge-Kutta method.

The numerical experiment, involving the Van der Pol equation with both well-prepared and unprepared initial conditions, demonstrated the robustness and efficiency of the proposed methods. The semi-implicit and implicit schemes showed excellent agreement with the reference solutions computed using MATLAB's ode15s solver. Detailed examinations at critical times revealed no significant differences between the results obtained from well-prepared and unprepared initial conditions, underscoring the AP (asymptotic preserving) nature of the schemes considered.

Furthermore, the adaptive time-step control proved crucial in managing the challenging boundary layers inherent in the Van der Pol problem, ensuring both accuracy and computational efficiency. The stability analyses confirmed the theoretical expectations, validating the use of these high-order Taylor series-based methods for practical applications.

Our findings suggest that the proposed schemes are not only reliable and accurate but also computationally efficient, making them suitable for a wide range of stiff ODE problems. Future work could explore the integration of more advanced time-step control techniques, such as the a posteriori Multi-dimensional Optimal Order Detection (MOOD) paradigm, to further enhance computational efficiency.

\section*{Acknowledgements} 
This research has received funding from the European Union’s NextGenerationUE – Project: Centro Nazionale HPC, Big Data e Quantum Computing, “Spoke 1” (No. CUP E63C22001000006). E. Macca was partially supported by GNCS No. CUP E53C23001670001 Research Project “Metodi numerici per le dinamiche incerte”. E. Macca and S. Boscarino would like to thank the Italian Ministry of Instruction, University and Research (MIUR) to support this research with funds coming from PRIN Project 2022  (2022KA3JBA, entitled “Advanced numerical methods for time dependent parametric partial differential equations and applications”). Sebastiano Boscarino has been supported for this work from Italian Ministerial grant PRIN 2022 PNRR “FIN4GEO: Forward and Inverse Numerical Modeling of hydrothermal systems in volcanic regions with application to geothermal energy exploitation.”, No. P2022BNB97. E. Macca and S. Boscarino are members of the INdAM Research group GNCS.

On behalf of all authors, the corresponding author states that there is no conflict of interest.

%

\bibliographystyle{plain}
\bibliography{biblio}

\end{document}